\documentclass[ieot]{birkjour}
\usepackage{pdfsync}
\usepackage{amsfonts}
\usepackage{enumerate} 
\usepackage{amsthm}
\usepackage{amsmath}
\usepackage{mathtools}
\usepackage{amscd}
\usepackage[active]{srcltx} % SRC Specials: DVI [Inverse] Search
\usepackage{latexsym}
\usepackage{amssymb}
\usepackage{enumitem}
\newtheorem{thm}{Theorem}[section]
\newtheorem{prop}[thm]{Proposition}
\newtheorem{cor}[thm]{Corollary}
\newtheorem*{cor*}{Corollary}
\newtheorem{lema}[thm]{Lemma}
\newtheorem*{lema*}{Lemma}

\numberwithin{equation}{section}
\theoremstyle{definition}
\newtheorem*{Def}{Definition}

\newenvironment{dem}{\vspace{1ex}\noindent{\it Proof.}\hspace{0.5em}}
{\hfill\qed\vspace{1ex}}

\newtheorem*{obs}{Remark}

\newtheorem*{obs*}{Remark}
\newtheorem*{thm*}{Theorem}
\newtheorem*{prop*}{Proposition}

\newtheoremstyle{dotless}{}{}{}{}{}{}{ }{}
\theoremstyle{dotless}

\newcommand{\PI}[2]{\left\langle \,#1 , #2\, \right\rangle}

\newcommand{\K}[2]{[ \,#1 , #2\, ]}

\newcommand{\vareps}{\varepsilon}

\newcommand{\St}{\mathcal{S}}

\newcommand{\HH}{\mathcal{H}}

\newcommand{\M}{\mathcal{M}}
\newcommand{\N}{\mathcal{N}}

\newcommand{\Q}{\mathcal{Q}}

\newcommand{\KK}{\mathcal{K}}
\newcommand{\mc}[1]{\mathcal{#1}}

\newcommand{\ol}{\overline}

\newcommand{\perpi}{[\perp]}

\begin{document}

%\begin{frontmatter}

\title[Operator LSP and Moore-Penrose inverses]{Operator least squares problems and Moore-Penrose inverses in Krein spaces}

%\author[Birkh\"auser]{Birkh\"{a}user Publishing Ltd.}
%
%\address{%
%	Viaduktstr. 42\\
%	P.O. Box 133\\
%	CH 4010 Basel\\
%	Switzerland}
%
%\email{info@birkhauser.ch}

\author[Contino]{Maximiliano Contino}
\address{Instituto Argentino de Matem\'atica ``Alberto P. Calder\'on''\\ Saavedra 15, Piso 3 (1083) Buenos Aires, Argentina\\ and \\ Departamento de Matem\'atica -- Facultad de Ingenier\'{\i}a -- Universidad de Buenos Aires\\ Paseo Col\'on 850 (1063) Buenos Aires, Argentina}
\email{mcontino@fi.uba.ar}

%\author[IAM]{Guillermina Fongi}
%\ead{gfongi@gmail.com}

\author[Maestripieri]{Alejandra Maestripieri}
\address{Instituto Argentino de Matem\'atica ``Alberto P. Calder\'on''\\ Saavedra 15, Piso 3 (1083) Buenos Aires, Argentina \\ and \\ Departamento de Matem\'atica -- Facultad de Ingenier\'{\i}a -- Universidad de Buenos Aires\\ Paseo Col\'on 850 (1063) Buenos Aires, Argentina}
\email{amaestri@fi.uba.ar}

\author[Marcantognini]{Stefania Marcantognini}
\address{Departamento de Matem\'atica -- Instituto Venezolano de Investigaciones Cient\'ificas \\ Km 11 Carretera Panamericana Caracas, Venezuela \\ and \\ Instituto Argentino de Matem\'atica ``Alberto P. Calder\'on''\\ Saavedra 15, Piso 3 (1083) Buenos Aires, Argentina}
\email{stefania.marcantognini@gmail.com}
%\cortext[ca]{Corresponding author.  \thedate}

\begin{abstract}
A Krein space $\HH$ and bounded linear operators $B, C$ on $\HH$ are given. Then, some min and max problems about the operators  \linebreak
$(BX - C)^{\#}(BX -C)$, where $X$ runs over the space of all bounded linear operators on 
$\HH$, are discussed. In each case, a complete answer to the problem, including solvability conditions and characterization of the solutions, is presented. Also, an adequate decomposition of $B$ is considered and the min-max problem is addressed. As a by-product the Moore-Penrose inverse of $B$ is characterized as the only solution of a variational problem. Other generalized inverses are described in a similar fashion as well.

\end{abstract}



%\end{frontmatter}

\maketitle

\section{Introduction}
Several least squares problems, especially in connection with the search of alternative $H^\infty$ algorithms in system and control theory, have been placed in the Krein space framework. Roughly speaking, the consideration of suitable space models for the set of observations data have brought into play indefinite metric spaces and, on the basis of the given information, least squares problems on those spaces. Some references from the nineties are \cite{HassibipartI, HassibipartII, HassibiIII}. Commonly those least squares estimations are formulated and solved in terms of vectors in Krein spaces and more often than not the vectors are set in the $(n+m)$-dimensional Minkowski space. 

We discuss, instead, least-squares problems for Krein space operators. The approach we opt for is taken from \cite{GiribetIndefinite, GiribetKrein}. Several arguments we present are adapted from \cite{Contino, Contino2}.

Our aim is dual: to study abstract least-squares-type problems for Krein space operators and to do so from a geometrical viewpoint. Pseudo-regularity plays a key role,  either as a technical tool to generalize some finite-dimensional indefinite metric space arguments to the general Krein space framework or as the natural  assumption to grant solvability. In that regard we should mention that every closed subspace of a Pontryagin space is pseudo-regular and, on the other hand, that every pseudo-regular subspace of a Krein space is the range of a normal projection. We may say so that pseudo-regular subspaces lie somewhere in between closed subspaces and regular subspaces. For more details on the subject see \cite{Gheondea} and \cite{Pseudo}.

Recall that, in the Hilbert space case, the Moore-Penrose inverse $B^\dagger$ of a given bounded linear operator $B$ satisfies the equations $BB^{\dagger}B=B$  and $B^{\dagger}BB^{\dagger}=B^{\dagger}$, and that both $BB^{\dagger}$ and $B^{\dagger}B$ are selfadjoint.  If $B$ is a closed range operator then $B^\dagger$ is known to be the unique minimal norm solution of a linear equation. In the Krein space framework the analysis of the existence of a generalized inverse  $B^\dagger$ such that $BB^{\dagger}$ and $B^{\dagger}B$ are selfadjoint -- with respect to the indefinite inner product -- was carried on by X. Mary~\cite{XavierMary}. He found out that  a bounded linear operator $B$ admits a unique bounded Moore-Penrose inverse if and only if both the range and null space of $B$ are regular subspaces. His treatment is exhaustive but it fails to include the variational characterization of the Moore-Penrose inverse.

Under the necessary and sufficient conditions given by Mary, 
we do identify $B^\dagger$ as the unique solution of a variational problem. Furthermore, when the projections associated to the generalized inverse are only required to be normal -- with respect to the indefinite inner products -- we prove that pseudo-regularity of the range $R(B)$ and null space $N(B)$ are necessary and sufficient conditions for a closed range $B$ to admit such a sort of generalized inverse. Matter-of-factly, if that is the case, there exists a whole family of generalized inverses which is in one to one correspondence with the set of  pairs $(Q, P)$ with $Q$ a normal projection onto $R(B)$ and $P$ a normal projection onto $N(B)$.  Besides, we characterize the generalized inverses in a variational way just as we do it for the Moore-Penrose inverse. 

The paper comprises five sections, six if this introductory section is included. Section 2 is a brief expository introduction to Krein spaces and operators on them and serves to fix the notation. It presents also some results that are needed in the following sections. In Section 3 the notion of indefinite inverse of an operator is defined, generalizing the concept of $W$-inverse introduced by Mitra and Rao for matrices in \cite{Mitra}, and extended later for operators acting on Hilbert spaces in \cite{WGI}. In Section 4 we deal with the problem of determining whether the 
$$\underset{X \in L(\HH)}{min} (BX-C)^{\#}(BX-C) \footnote{$L(\HH)$ stands for the space of all the bounded linear operators from $\HH$ to $\HH$.}$$ 
exists for $B$ a given closed range bounded operator and $C$ either the identity operator or any given bounded operator. The solutions to these problems are characterized as the indefinite inverses of $B.$ The results about this indefinite minimization problem and their counterparts for the symmetric maximization problem are applied in Section 5 where $B$ is factorized as $B = B_+ + B_-$ and the min-max problem
$$\underset{X \in L(\HH)}{min} \underset{Y \in L(\HH)}{max}  (B_+X+B_-Y-C)^{\#}(B_+X+B_-Y-C)$$
is addressed. Section 6 contains the main results about the Moore-Penrose inverse and the generalized inverses of a given Krein space operator.

\section{Preliminaries}

In the following all Hilbert spaces are complex and separable. If $\HH$ and $\KK$ are Hilbert spaces, $L(\HH, \KK)$ stands for the space of all the bounded linear operators from $\HH$ to $\KK$ and $CR(\HH, \KK)$ for the subset of $L(\HH, \KK)$ comprising all the operators with closed ranges. When $\HH = \KK$ we write, for short, $L(\HH)$  and $CR(\HH)$. 
%Besides, we write $GL(\HH)$ to indicate the group of all the invertible operators in the algebra $L(\HH).$ 
The range and null space of any given $A \in L(\HH, \KK)$ are denoted by $R(A)$ and $N(A)$, respectively. 

The direct sum of two closed subspaces $\M$ and $\N$ of $\HH$ is represented by $\M \dot{+} \N.$ 
If $\HH$ is decomposed as $\HH=\M \dot{+} \N,$ the projection onto $\M$ with null space $\N$ is denoted $P_{\M {\mathbin{\!/\mkern-3mu/\!}} \N}$ and abbreviated $P_{\M}$ when $\N = \M^{\perp}.$ In general, $\Q$ is used to indicate  the subset of all the oblique projections in  $L(\HH),$ namely, $\Q:=\{Q \in L(\HH): Q^{2}=Q\}.$

The following  is a well-known result about range inclusion and factorizations of operators. We will refer to it along the paper.
\begin{thm} [Douglas' Theorem \cite{Douglas}] \label{teo Douglas} Let $Y, Z \in L(\HH)$. Then $R(Z)\subseteq R(Y)$ if and only if there exists $D\in L(\HH)$ such that $Z=YD.$
\end{thm}

\subsection*{\textbf{Krein Spaces}}

Although familiarity with operator theory on Krein spaces is presumed, we hereafter include some basic notions. Standard references on Krein spaces and operators on them are \cite{AndoLibro, Azizov, Bognar}. We also refer to \cite{DR, DR1} as authoritative accounts of the subject.

Consider a linear space $\HH$ with an indefinite metric, i.e., a sesquilinear Hermitian form $\K{ \ }{ \ }.$  A vector $x \in \HH$ is said to be {\sl{positive}} if $\K{x}{x} > 0.$ A subspace $\St$ of $\HH$ is {\sl{positive}} if every $x \in \St,$ $x \not =0,$ is a positive vector. {\sl{Negative, nonnegative, nonpositive}} and {\sl{neutral}} vectors and subspaces are defined likewise. 

We say that two closed subspaces $\M$ and $\N$ are {\sl{orthogonal}}, and we write $\M  \ \perpi \ \N,$ if $\K{m}{n}=0 \mbox{ for every } m \in \M \mbox{ and } n \in \N.$ We denote the orthogonal direct sum of two closed subspaces $\M$ and $\N$ by $\M \ [\dotplus] \ \N.$  

Given any subspace $\St$ of $\HH,$ the {\sl{orthogonal companion}} of $\St$ in $\HH$, say $\St^{\perpi}$, is defined as 
$$\St^{\perpi}:=\{ x \in \HH: \K{x}{s}=0 \mbox{ for every } s \in \St \}.$$ The isotropic part $\St^o := \St \cap \St^{\perpi}$ can be a non-trivial subspace. 

An indefinite metric space $(\HH, \K{ \ }{ \ })$ is a {\sl{Krein space}} if $\HH$ admits a decomposition into an orthogonal direct sum in the form
\begin{equation} \label{fundamentaldecom}
\HH=\HH_+ \ [\dotplus] \ \HH_-
\end{equation} 
where $(\HH_+, \K{ \ }{ \ })$ and $(\HH_-, -\K{ \ }{ \ })$ are Hilbert spaces. Any decomposition with these properties is called a {\sl{fundamental decomposition}} of  $\HH.$ 
%Notice that any Hilbert space is itself a Krein space.

Given a Krein space $(\HH, \K{ \ }{ \ })$ with a fundamental descomposition $\HH=\HH_+ \ [\dotplus] \ \HH_-,$ the (orthogonal) direct sum of the Hilbert spaces $(\HH_+, \K{ \ }{ \ })$ and $(\HH_-, -\K{ \ }{ \ })$ is a Hilbert space. It is  denoted by $(\HH, \PI{ \ }{ \ }).$ Notice that the inner product $\PI{ \ }{ \ }$ and the corresponding quadratic norm $\Vert \ \Vert$ depend on the fundamental decomposition. 
A subspace $\St$ of $\HH$ is called {\sl{uniformly positive}} if, for some Hilbert space inner product $\PI{ \ }{ \ }$ on $\HH$, there exists $\vareps >0$ such that $\K{s}{s} \geq \vareps \Vert s \Vert^2 \mbox{ for every } s \in \St.$  {\sl{Uniformly negative}} subspaces are defined in a similar fashion.

Every fundamental decomposition of $\HH$ has an associated {\sl{signature operator}}, to wit, $J:=P_{+} - P_{-}$ where $P_{\pm}:=P_{\HH_{\pm}  {\mathbin{\!/\mkern-3mu/\!}} \HH_{\mp}}.$ The indefinite metric and the inner product corresponding to a fundamental decomposition of $\HH$ with signature operator $J$ are related to each other by 
$$\PI{x}{y}=\K{Jx}{y} \mbox{ for every } x, y \in \HH.$$ 

If $\HH$ is a Krein space, $L(\HH)$ stands for the vector space of  all the linear operators on $\HH$ which are bounded in an associated Hilbert space $(\HH, \PI{ \ }{ \ }).$ Since the norms generated by different fundamental decompositions of a Krein space $\HH$ are equivalent, see, for instance, \cite[Theorem 7.19]{Azizov}, it comes that $L(\HH)$ does not depend on the chosen underlying Hilbert space.

Given $T \in L(\HH),$ $T^{\#}$ is the unique operator satisfying $$\K{Tx}{y}=\K{x}{T^{\#}y} \mbox{ for every } x, y \in \HH.$$
An operator $T \in L(\HH)$ is said to be {\sl{selfadjoint}} if $T=T^{\#}.$ 

A {\sl{positive}} operator $T\in L(\HH)$ satisfies $\K{Tx}{x} \geq 0 \mbox{ for every } x \in \HH.$ The notation 
$S \leq T$ signifies that $T-S$ is positive.

A (closed) subspace $\St$ of a Krein space $\HH$ is a {\sl{regular subspace}} if $\HH=\St \ [\dotplus] \ \St^{\perpi}.$ Equivalently, $\St$ is a regular subspace if it is the range of a selfadjoint projection, i.e., there exists $Q \in \Q$ such that $Q=Q^{\#}$ and $R(Q)=\St$ (see \cite[Proposition 1.4.19]{Azizov}).

In \cite[Theorem 2.3]{Ando}, T.~Ando proved that any selfadjoint projection on a Krein space can be decomposed as the sum of two selfadjoint projections with uniformly definite ranges. See also \cite{Hassi, KreinSzeged}.

\begin{thm} \label{THM2.3} Let $(\HH, \K{ \ }{ \ })$ be a Krein space and let $Q$ be a selfadjoint projection. Then $Q$ can be written as $$Q=Q_+ +Q_-,$$ where $Q_+$ and $Q_-$ are  selfadjoint projections such that $R(Q_+)$ is uniformly positive, $R(Q_-)$ is uniformly negative and $Q_+Q_-=Q_-Q_+=0.$ Moreover, each fundamental decomposition of $\HH$  provides a (unique) decomposition of $Q$ in such a manner.
\end{thm}

%\begin{Def}  Let $(\HH, \K{ \ }{ \ })$ be a Krein space  and consider a fundamental decomposition $\HH=\HH_+ \ [\dotplus] \ \HH_-$ and $\St \subseteq \HH$ a closed subspace. If $J= P_+ - P_-$ is the associated fundamental symmetry, the Gram operator of $\St$ with respect to $J$ is 
%	$$G_{\St}=P_{\St}J|_{\St},$$
%	
%where $P_{\St}$ is the orthogonal projection onto $\St$ in the associated Hilbert space $(\HH, \PI{ \ }{ \ }).$	
%	
%For every $x, y \in \St$ we have that $\K{x}{y}=\PI{G_{\St}x}{y}.$
%\end{Def}

The next lemma shows that every closed subspace of a Krein space can be decomposed as the orthogonal direct sum of a closed positive subspace and a closed nonpositive subspace (see \cite[Theorem 6.4]{Azizov}, \cite[Chapter V, Theorem 3.1]{Bognar}).

\begin{lema} \label{lemaDecom} Let $(\HH, \K{ \ }{ \ })$ be a Krein space  with fundamental descomposition $\HH=\HH_+ \ [\dotplus] \ \HH_-$ and corresponding Hilbert space inner product $\PI{ \ }{ \ }$. Let $\St$ be a closed subspace of $\HH$. Then $\St$ can be represented uniquely as the orthogonal direct sum of a closed positive subspace $\St_+$ and a closed nonpositive subspace $\St_-$, i.e.,
$$\St = \St_+ \ [\dotplus] \ \St_-.$$
Furthermore, $\langle\St_+,\St_-\rangle = \{0\}$. 
\end{lema}

%\begin{dem} Let $G:=G_{\St}$ be the Gram operator of $\St.$ Then, since $G$ is selfadjoint in the associated Hilbert space $(\HH, \PI{ \ }{ \ }),$ there exist operators $G_1, G_2 \in L(\HH)^+,$ such that $G=G_1-G_2,$  $R(G_1) \perp R(G_2),$
%and $\ol{R(G)}=\ol{R(G_1)} \oplus \ol{R(G_2)}.$ Therefore
%$$\St = N(G) \oplus \ol{R(G)} = N(G) \oplus \ol{R(G_1)} \oplus \ol{R(G_2)}.$$
%Observe that $\St^{o}=N(G).$ In fact, $x \in N(G)$ if and only if $x \in \St$ and $G=0,$ if and only if $x \in \St$ and $\PI{Gx}{y}=0, \mbox{ for every } y \in \St,$ if and only if $x \in \St$ and $\K{x}{y}=0, \mbox{ for every } y \in \St,$ if and only if $x \in \St^{o}.$
%Let $\St_1=\ol{R(G_1)}$ and $\St_2=\ol{R(G_2)},$ then $\St_1$ is $J$-positive and $\St_2$ is $J$-negative.  In fact, let $x \in \St_1 \setminus \{0\},$ then
%$$\K{x}{x}=\PI{Gx}{x}=\PI{GP_{\ol{R(G_1)}}x}{x}=\PI{G_1x}{x} \geq 0,$$ moreover $\PI{G_1x}{x} >0,$ because $\PI{G_1x}{x}=0$ if and only if $\Vert (G_1)^{1/2} x \Vert=0,$ if and only if $x \in N(G_1) \cap R(G_1) =\{0\},$ therefore $\St_1$ is $J$-positive. Analogously, it can be proved that $\St_2$ is $J$-negative.
%
%Finally, observe that $\St_1 \ \perpi \ \St_2,$ in fact, let $x \in \St_1$ and $y \in \St_2$ then
%$$\K{x}{y}=\PI{Gx}{y}=0.$$
%
%Now, define $\St_+=\St_1,$ and $\St_-=\St^{0} \ [\dotplus] \ \St_2,$ then $\St_+$ and $\St_-$ are closed $J$-positive and $J$-nonpositive subspaces respectively, $S_+ \perp \St_-$ and
%
%$$\St=\St^{o} \ [ \dotplus ] \ \St_1 \ [ \dotplus ] \ \St_2 = \St_+ \ [\dotplus ] \ \St_-.$$
%\end{dem}
% 
 
In \cite{GiribetIndefinite, GiribetKrein} least squares problems in the indefinite metric setting were studied.  From those references we recall the definition of indefinite least squares solution. 

\begin{Def} Let $(\HH, \K{ \ }{ \ })$ be a Krein space and let $B\in CR(\HH).$ We say that $u \in \HH$ is an {\sl{indefinite least squares solution}} (ILSS) of $Bz=x$ if
	$$\K{x-Bu}{x-Bu}\leq \K{x-Bz}{x-Bz} \ \mbox{ for every } x, z \in \HH.$$
\end{Def}

We conclude this section by stating
necessary and sufficient conditions for the existence of an ILSS of the equation $Bz=x$. We refer to \cite[Chapter I, Theorem 8.4]{Bognar} where a proof of the result is given.

\begin{lema} \label{Lemma3.1} Let $(\HH, \K{ \ }{ \ })$ be a Krein space and let $B\in CR(\HH).$ Then $u\in \HH$ is an ILSS of the equation $Bz=x$ if and only if $R(B)$ is nonnegative and $x-Bu \in R(B)^{\perpi}.$
\end{lema}

 \section{Indefinite inverses in Krein spaces}
 
In \cite{Mitra} S. K. Mitra and C. R. Rao introduced the notion of the $W$-inverse of a matrix for a given positive weight $W.$ Later, in \cite{WGI} and \cite{Contino}, the concept was extended to Hilbert space operators, specifically, given a Hilbert space $(\HH, \PI{ \ }{ \ })$, a positive operator $W\in L(\HH)$ 
%(i.e., such that $\PI{Wx}{x} \geq 0 \mbox{ for every } x \in \HH$), 
and an operator $B \in CR(\HH),$ 
a {\sl{$W$-inverse}} of $B$ is defined  to be an operator $X_0 \in L(\HH)$ such that, for each $x \in \HH,$ $X_0x$ is a weighted least squares solution of $Bz=x,$ i.e., so that
$$\PI{W(BX_0x-x)}{BX_0x-x}\leq\PI{W(Bz-x)}{Bz-x} \mbox{ for every } z \in \HH.$$  

In \cite{WGI} it was proved that $X_0$ is a W-inverse of $B$ if and only if \linebreak $B^*W(BX_0-I)=0$ or, equivalently, $X_0$ satisfies the identities $W(BX_0)^2=WBX_0=(BX_0)^*W.$ The first equality means that $BX_0$ is a {\sl{$W$-projection}} while the second says that $BX_0$ is {\sl{$W$-selfadjoint}}, see \cite{WGI}.
 
We extend the definition to Krein spaces in the following way.
From now on, $(\HH, \K{ \ }{ \ })$ stands for  a Krein space.

\begin{Def} Let $B\in CR(\HH).$ An operator $X_0 \in L(\HH)$ is an {\sl{indefinite inverse}} of $B$ if $X_0$ is a solution of $$B^{\#} (BX-I)=0.$$
\end{Def}

\begin{prop} \label{PropJinvReg} Let $B\in CR(\HH).$  Then $B$ admits an indefinite inverse if and only if $R(B)$ is regular.
\end{prop} 
\begin{dem} Suppose that $X_0$ is an indefinite inverse of $B$ so that $B^{\#}(BX_0-I)=0.$ Then, for every $x \in \HH,$ $BX_0x-x \in N(B^{\#})=R(B)^{\perpi}$ and, therefore, $x \in R(B) + R(B)^{\perpi}$. Whence $\HH=R(B) + R(B)^{\perpi}.$ As, also, $\{0\}=R(B) \cap R(B)^{\perpi}$, it comes that $\HH=R(B) \ [\dotplus] \ R(B)^{\perpi}$ or, accordingly, that $R(B)$ is regular.
	
Conversely, if $R(B)$ is regular then $\HH=R(B) [ \dotplus ] R(B)^{\perpi}.$ So, by applying $B^{\#}$, it results that  
$
R(B^{\#}) = R(B^{\#}B).
$
From here and by Douglas' Theorem (Theorem \ref{teo Douglas}), it follows that the equation $B^{\#} (BX-I)=0$ admits a solution or, equivalently, that $B$ has  an indefinite inverse.
\end{dem}

%\vspace{0,3cm} 
It results from the proof of Proposition \ref{PropJinvReg} that, for any $B \in CR(\HH),$ $R(B)$ is regular if, and only if, 
$R(B^{\#}) = R(B^{\#}B).$
In this case,
$N(B)=N(B^{\#}B).$
%See also \cite{XavierMary}.

The next proposition characterizes the indefinite inverses of $B \in L(\HH)$ when $R(B)$ is regular.

\begin{prop} \label{PropJ1b} Let $B\in L(\HH)$. Assume that $R(B)$ is regular. Then the following conditions are equivalent:
	\begin{enumerate}
		\item [i)] $X_0$ is an indefinite inverse of $B$, 
		\item [ii)] $X_0$ is a solution of the equation $BX=Q,$ where $Q$ is the selfadjoint projection onto $R(B),$
		\item [iii)] $X_0$ is an inner inverse of $B,$ i.e., $BX_0B=B,$ and $(BX_0)^{\#}=BX_0.$
 \end{enumerate}
 Moreover, if $R(B)$ is also uniformly positive, conditions $i), ii), iii)$ are also equivalent to:
 	\begin{enumerate}
	\item [iv)] For every $x \in \HH,$  $X_0x$ is an ILSS of $Bz=x.$ 
	 \end{enumerate}
	A similar statement holds if $R(B)$ is uniformly negative.

\end{prop}
\begin{dem} \mbox{$i) \Rightarrow ii):$} Notice, first, that $B^{\#}=B^{\#}Q,$ since $B=QB$ and $Q=Q^{\#}.$  Whence $B^{\#}(BX_0-I)=0$ implies $B^{\#}(BX_0-Q)=0.$ Therefore, $R(BX_0-Q) \subseteq N(B^{\#}) \cap R(Q) = N(Q) \cap R(Q) = \{0\}.$ Thus $BX_0=Q.$
	
	\mbox{$ii) \Rightarrow iii):$} If $BX_0=Q$ then $BX_0B=QB=B$ and $(BX_0)^{\#}=Q^{\#}=Q=BX_0.$
	
	\mbox{$iii) \Rightarrow i):$}  Suppose that $BX_0B=B$ and $(BX_0)^{\#}=BX_0.$ Then \linebreak $B^{\#}(BX_0-I)=B^{\#}(X_0^{\#}B^{\#}-I)=B^{\#}X_0^{\#}B^{\#}-B^{\#}=0.$
	
	$i) \Leftrightarrow iv):$  $X_0$ is an indefinite inverse of $B$ if, and only if, $B^{\#}(BX_0-I)=0$ if, and only if, for every $x \in \HH,$  $BX_0x-x \in R(B)^{\perpi}.$ Since $R(B)$ is nonnegative, Lemma \ref{Lemma3.1} gives the equivalence.	
\end{dem}

%\vspace{0,3cm} 
We point out that, when $R(B)$ is closed and uniformly definite, the original definition of $W$-inverse for the indefinite metric is retrieved directly from item $iv)$ of the last proposition. 

The more general concept of the {\sl{indefinite inverse}} of $B$ {\sl{in the range of}} $C$ is given next.

\begin{Def} Let  $B\in CR(\HH)$ and $C \in L(\HH).$ An operator $X_0 \in L(\HH)$ is an {\sl{indefinite inverse}} of $B$ {\sl{in}} $R(C)$ if $X_0$ is a solution of $$B^{\#} (BX-C)=0.$$
\end{Def}

\begin{prop} \label{PropJinvReg2} Let $B\in CR(\HH)$ and $C \in L(\HH).$  $B$ has an indefinite inverse in $R(C)$ if and only if $R(C) \subseteq R(B) + R(B)^{\perpi}.$ 
\end{prop} 
\begin{dem} Suppose that $X_0$ is an indefinite inverse of $B$ in $R(C).$ Then \linebreak $B^{\#}(BX_0-C)=0.$ So, if $x \in \HH$ then $BX_0x-Cx \in N(B^{\#})=R(B)^{\perpi}$ and, therefore, $Cx \in R(B) + R(B)^{\perpi}$. Thus, $R(C) \subseteq R(B) + R(B)^{\perpi}$  and the result follows.
	
Conversely, if $R(C) \subseteq R(B) + R(B)^{\perpi}$ then
$$R(B^{\#}C) \subseteq R(B^{\#}B).$$ 
Here, as before in the proof of Proposition \ref{PropJinvReg}, Douglas' Theorem is applied to grant that the equation $B^{\#} (BX-C)=0$ admits a solution or, equivalently, that $B$ has an indefinite inverse in $R(C).$
\end{dem}

%\vspace{0,3cm} 

\begin{cor} \label{CorJinvC} Let $B \in CR(\HH)$ and $C \in L(\HH).$ 
%such that $R(C) \subseteq R(B) + R(B)^{\perpi}.$  
%Then $X_0$ is an indefinite inverse of $B$ in $R(C),$ if and only if  $$X_0=D+T,$$ with $D$ a solution of the normal equation $B^{\#}(BX-C)=0$ and $R(T) \subseteq N(B^{\#}B).$
If $R(B)$ is regular then $X_0$ is an indefinite inverse of $B$ in $R(C)$ if and only if $X_0$ is a solution of the equation $BX=QC,$ with $Q$ the selfadjoint projection onto $R(B).$
\end{cor}
\begin{dem} Suppose that $R(B)$ is regular. Then $R(B^{\#}) = R(B^{\#}B)$ or, equivalently, $N(B)=N(B^{\#}B).$  

If $B^{\#}(BX_0-C)=0$ then $B^{\#}(BX_0-QC)=0,$ for $B=QB$, $Q=Q^{\#}$ and, consequently, 
$B^{\#}C=B^{\#}QC.$ Hence $R(BX_0-QC) \subseteq N(B^{\#}) \cap R(Q) = N(Q) \cap R(Q) = \{0\}.$ So $BX_0=QC.$ 
\end{dem}

From the last corollary we have that, when $R(B)$ is regular, the set of indefinite inverses of $B$ in $R(C)$ is the affine manifold 
$$X_0+L(\HH, N(B)),$$ 
with $X_0$ any solution of the equation $BX=QC.$

\begin{prop} \label{PropJ2} Let $B\in CR(\HH)$ and $C \in L(\HH)$ satisfy that $R(B)$ is nonnegative and $R(C) \subseteq R(B) + R(B)^{\perpi}.$ Then $X_0$ is an indefinite inverse of $B$ in $R(C)$ if and only if, for every $x \in \HH,$ $X_0x$ is an ILSS of $Bz=Cx,$ i.e., 
	$$\K{Cx-BX_0x}{Cx-BX_0x}\leq \K{Cx-Bz}{Cx-Bz} \ \mbox{ for every } z \in \HH.$$
\end{prop}
\begin{dem} $X_0$ is an indefinite inverse of $B$ in $R(C)$ if, and only if, $B^{\#}(BX_0-C)=0$ if, and only if,   $BX_0x-Cx \in R(B)^{\perpi}$ for every $x \in \HH.$ Since $R(B)$ is supposed to be nonnegative, Lemma \ref{Lemma3.1} can be applied to get the equivalence.
\end{dem}

\section{Indefinite least squares problems}
To state the next problems let us recall that the order is the one induced by the positive operators in $(\HH, \K{ \ }{ \ })$: given two operators $S, T \in L(\HH),$ $S \leq T$ whenever $T-S$ is positive. 

Consider the following problem: given 
two operators $B\in CR(\HH)$ and $C \in L(\HH),$ determine the existence of 
%\begin{equation}
%\underset{X \in L(\HH)}{min_{\leq_{J}}} (BX-I)^{\#}(BX-I). \label{eq41}
%\end{equation}
%
%\vspace{0,3cm} 
%More generally, we study the following problem:
%let $(\HH, \K{ \ }{ \ })$ be a Krein space, $B\in CR(\HH)$ and $C \in L(\HH),$ analize the existence of
\begin{equation}
\underset{X \in L(\HH)}{min} (BX-C)^{\#}(BX-C). \label{eq51}
\end{equation}

\begin{Def} Let $B\in CR(\HH)$ and $C \in L(\HH).$ We say that $X_0 \in L(\HH)$ is an {\sl{indefinite minimum solution}} (ImS) of $BX-C=0$ if 	
\begin{equation} \label{eqMin2}
(BX_0-C)^{\#}(BX_0-C)=\underset{X \in L(\HH)}{min} (BX-C)^{\#}(BX-C).
\end{equation}
\end{Def}

In a similar fashion, the analogous maximization problem can be considered.
From now on, we only address the problem related to the existence of \eqref{eq51}. The arguments we present in dealing with problem  \eqref{eq51} can be adapted to the maximum problem. In particular, each of the ``min" results we include in this section can be easily modified to get its   ``max" counterpart.

\begin{thm}  \label{Teo1011} Let $B\in CR(\HH)$ and $C \in L(\HH).$ 
Then there exists an ImS of $BX-C=0$ if and only if  $R(C) \subseteq R(B) + R(B)^{\perpi}$ and $R(B)$ is nonnegative.
	
\end{thm}
\begin{dem} Suppose that $X_0\in L(\HH)$ is an ImS of $BX-C=0,$ so that 
	$$ \K{(BX_0-C)x}{(BX_0-C)x} \leq \K{(BX-C)x}{(BX-C)x}$$
	$ \mbox{for every } x \in \HH \mbox{ and every } X \in L(\HH).$
	Let $z \in \HH$ be arbitrary. Then, for every $x \in \HH \setminus \{0\},$ there exists $X \in L(\HH)$ such that $z=Xx.$ Therefore
	$$ \K{(BX_0-C)x}{(BX_0-C)x} \leq \K{Bz-Cx}{Bz-Cx} \mbox{ for every } x, z \in \HH.$$
	So, for every $x \in \HH,$ $X_0x$ is an ILSS of $Bz=Cx.$ By Lemma \ref{Lemma3.1}, we get that $R(C) \subseteq R(B) + \ R(B)^{\perpi}$ and $R(B)$ is nonnegative. Furthemore, by Proposition \ref{PropJ2},  we have that $X_0$ is an indefinite inverse of $B$ in $R(C).$
	
	Conversely, if $R(C) \subseteq R(B) + R(B)^{\perpi}$ and $R(B)$ is nonnegative then,  by Proposition \ref{PropJinvReg2}, $B$ admits an indefinite inverse in $R(C).$ Now, if $X_0$ is an indefinite inverse of $B$ in $R(C)$ then, by Proposition \ref{PropJ2},
	$$\K{(BX_0-C)x}{(BX_0-C)x} \leq \K{Bz-Cx}{Bz-Cx} \mbox{ for every } x, z \in \HH.$$ 
	Given $x \in \HH$ and  $X \in L(\HH),$ set $z=Xx,$ so that 
	$$\K{(BX_0-C)x}{(BX_0-C)x} \leq \K{(BX-C)x}{(BX-C)x}$$
	$ \mbox{for every } x \in \HH  \mbox{ and every } X \in L(\HH).$
	Now it becomes clear that $X_0$ is an ImS of the equation $BX-C=0,$ as required to complete the proof.
\end{dem}

%\bigskip

%	If $R(C) \subseteq R(B) + R(B)^{\perpi},$ $R(B)$ is nonnegative and $X_0$ is a $J$-inverse of $B$ in $R(C)$, by Corollary \ref{CorJinvC}, $$X_0=(B^{\#}B)^{\dagger}B^{\#}C+T,$$ with $R(T) \subseteq N(B^{\#}B),$ then 
%	\begin{equation} \label{minimoC}
%	\underset{X \in L(\HH)}{min_{\leq_{J}}} (BX-C)^{\#}(BX-C)= (B(B^{\#}B)^{\dagger}B^{\#}C-C)^{\#}(B(B^{\#}B)^{\dagger}B^{\#}C-C),
%	\end{equation}
%	where we used $C^{\#}BT=0,$ because $R(C) \subseteq R(B) + R(B)^{\perpi}.$

%\bigskip
In the proof of Theorem \ref{Teo1011} the $X_0$'s in \eqref{eqMin2}, that is, the solutions of the problem related to \eqref{eq51}, were characterized. Indeed:

\begin{cor} \label{Prop103} Let $B\in CR(\HH)$ and $C \in L(\HH)$ satisfy that $R(B)$ is nonnegative and $R(C) \subseteq R(B) + R(B)^{\perpi}.$ Then $X_0$ is an ImS of $BX-C=0$ 
%a solution of Problem \eqref{eq51}%
 if and only if $X_0$ is an indefinite inverse of $B$ in $R(C),$ i.e., $X_0$ is solution of the normal equation $$B^{\#}(BX-C)=0.$$

\end{cor}

%\bigskip
From the last result we have that the set of indefinite inverses of ImS of $BX-C=0$ is the affine manifold 
$$X_0+L(\HH, N(B^{\#}B)),$$ 
with $X_0$ any indefinite inverse of $B$ in $R(C).$

The next two corollaries follow from Theorem \ref{Teo1011} as well.

%\vspace{0,3cm} 
\begin{cor} \label{Cor1012} Let $B\in CR(\HH).$ Then there exists an ImS of $BX-C=0$ for every $C \in L(\HH)$ if and only if $R(B)$ is uniformly positive.
	
	In this case,
	$$\underset{X \in L(\HH)}{min} (BX-C)^{\#}(BX-C)=C^{\#}(I-Q)C$$ 
	where $Q$ is the selfadjoint projection onto $R(B).$
\end{cor}
\begin{dem} Assume that, for every $C \in L(\HH),$  there exists an ImS of $BX-C=0.$
In particular, there exists an ImS of $BX-I=0.$ Then, by Theorem \ref{Teo1011}, $R(B)$ is regular and nonnegative. Whence $R(B)$ is uniformly positive.
	
Conversely, if $R(B)$ is uniformly positive then, for every $C \in L(\HH),$ we have that $R(C) \subseteq \HH=R(B)+R(B)^{\perpi}.$ Hence, by Theorem \ref{Teo1011}, for every $C \in L(\HH),$ there exists  an ImS $X_0\in L(\HH)$ of $BX-C=0$ or, equivalently,
$X_0$ is a solution of the normal equation $B^{\#}(BX-C)=0$ (see Corollary \ref{Prop103}).
In this case, since $R(B)$ is regular, Corollary \ref{CorJinvC} gives that $BX_0=QC.$ Therefore,
$$\underset{X \in L(\HH)}{min} (BX-C)^{\#}(BX-C)=(BX_0-C)^{\#}(BX_0-C)=C^{\#}(I-Q)C.$$
\end{dem}

\begin{cor}  \label{Teo101} Let $B\in CR(\HH).$ Then there exists an ImS of  $BX-I=0$ if and only if $R(B)$ is uniformly positive.
	
	In this case, the ImS of  $BX-I=0$ are the indefinite inverses of $B$ and
	$$\underset{X \in L(\HH)}{min} (BX-I)^{\#}(BX-I)= I-Q$$ 
	where $Q$ is the selfadjoint projection onto $R(B).$
\end{cor}

%\bigskip
\begin{obs}  By mimicking the arguments in the proof of Theorem \ref{Teo1011}, a similar result can be proved for operators acting between different Krein spaces. More precisely, let $(\HH, \K{ \ }{ \ }_{\HH}),$ $(\KK, \K{ \ }{ \ }_{\KK})$ and $(\mc{F}, \K{ \ }{ \ }_{\mc{F}})$ be  Krein spaces. Let $B\in CR(\HH,\KK)$ and $C \in L(\mc{F},\KK).$ Then there exists $X_0 \in L(\mc{F},\HH)$ such that
	$$\underset{X \in L(\mc{F},\HH)}{min} (BX-C)^{\#}(BX-C)=(BX_0-C)^{\#}(BX_0-C)$$
	 if and only if  $R(C) \subseteq R(B) + R(B)^{{\perpi}_{\KK}}$ and $R(B)$ is nonnegative.
	%In fact, there exists $X_0 \in L(\HH)$ such that
	%$\underset{X \in L(\HH)}{min_{\leq_{J}}} (BX-C)^{\#}(BX-C)=(BX_0-C)^{\#}(BX_0-C),$ if and only if,
	%$$ \K{(BX_0-C)x}{(BX_0-C)x}_{\KK} \leq \K{(BX-C)x}{(BX-C)x}_{\KK}, \mbox{ for every } x \in \HH \mbox{ and every } X \in L(\HH).$$ 
	%
	%Let $z \in \HH$ be arbitrary, then for every $x \in \HH \setminus \{0\},$ there exists $X \in L(\HH)$ such that $z=Xx,$ therefore
	%$$ \K{(BX_0-C)x}{(BX_0-C)x}_{\KK} \leq \K{Bz-Cx}{Bz-Cx}_{\KK}, \mbox{ for every } x, z \in \HH,$$
	%then, for every $x \in \HH,$ $X_0x$ is an ILSS of $Bz=Cx,$ then by Lemma \ref{Lemma3.1}, $R(C) \subseteq R(B) + \ R(B)^{\perpi}$ and $R(B)$ is $J_{\KK}$-nonnegative. 
	%
	%Conversely, if $R(C) \subseteq R(B) + R(B)^{\perpi}$ and $R(B)$ is $J_{\KK}$-nonnegative, applying $B^{\#} \in L(\KK,\HH)$ to both sides of the inclusion, we have
	%$$R(B^{\#}C) \subseteq R(B^{\#}B)$$ and by Theorem \ref{teo Douglas}, the equation $B^{\#} (BX-C)=0,$ admits a solution, then $B$ admits a $J$-inverse in $R(C).$ Let $X_0$ be a $J$-inverse of $B$ in $R(C),$ then by Proposition \ref{PropJ2},
	%$$\K{(C-BX_0)x}{(C-BX_0)x}_{\KK} \leq \K{Cx-Bz}{Cx-Bz}_{\KK}, \mbox{ for every } x, z \in \HH,$$ in particular, given $x \in \HH$ and  $X \in L(\HH),$ consider $z=Xx.$ Then 
	%$$\K{(C-BX_0)x}{(C-BX_0)x}_{\KK} \leq \K{(C-BX)x}{(C-BX)x}_{\KK}, \mbox{ for every } x \in \HH  \mbox{ and every } X \in L(\HH),$$
	%and $X_0 \in L(\HH)$ is such that
	%$\underset{X \in L(\HH)}{min_{\leq_{J}}} (BX-C)^{\#}(BX-C)=(BX_0-C)^{\#}(BX_0-C).$
\end{obs}

\subsection{Indefinite least squares problems: the pseudo-regular case}

A (closed) subspace $\St$ of a Krein space $\HH$ is called a {\sl{pseudo-regular subspace}} if the algebraic sum $\St + \St^{\perpi}$ is closed. Observe that, this is equivalent to the equality $(\St^0)^{\perpi}=\St + \St^{\perpi},$ see \cite{Gheondea}. Also, $\St$ is a pseudo-regular subspace if and only if $\St$ is the range of a normal projection, i.e., there exists $Q \in \Q$ such that $QQ^{\#}=Q^{\#}Q$ and $R(Q)=\St$ (see \cite[Theorem 4.3]{Pseudo}).
Unlike normal projections in Hilbert spaces, a normal projection in a Krein space need not be selfadjoint.
In what follows $\Q_{\St}$ stands for the set of normal projections onto the pseudo-regular subspace $\St,$ i.e.,
 $$\Q_{\St}:=\{Q \in L(\HH): Q^{2}=Q, \ QQ^{\#}=Q^{\#}Q, \ R(Q)=\St\}.$$ 
The set $\Q_{\St}$ has infinite elements, unless $\St$ is regular. See \cite{Pseudo} for further details on the subject.

Let $B \in CR(\HH),$ the next results relate the pseudo-regularity of $R(B)$ to the indefinite inverse of $B$ in $R(C)$ and the $ImS$ of  $BX-C=0.$

The next lemma, stated in \cite[Remark 2.1]{GiribetIndefinite}, will be useful when dealing with pseudo-regular ranges.

\begin{lema} \label{remark}
 Let $\St$ be a pseudo-regular subspace of $\HH.$
If $Q \in \Q_{\St}$ then $$Q^{\#}(I-Q)y=0 \mbox{ if and only if } y \in \St + \St^{\perpi}.$$
\end{lema}  

\begin{dem}  Let $Q \in \Q_{\St}.$ If $y \in \St + \St^{\perpi},$ by \cite[Remark 2.1]{GiribetIndefinite}, we have that $Q^{\#}(I-Q)y=0.$

Conversely, if $Q^{\#}(I-Q)y=0,$ since $Q^{\#}(I-Q) \in \Q,$ $y \in N(Q^{\#}(I-Q))=\St + \St^{\perpi}.$
\end{dem}

\begin{prop} \label{PropJinvPseudo} Let $B \in CR(\HH)$ and $C \in L(\HH).$ 
If $R(B)$ is pseudo-regular and $R(C) \subseteq R(B) + R(B)^{\perpi},$ then $X_0$ is an indefinite inverse of $B$ in $R(C)$ if and only if $R(BX_0-QC) \subseteq R(B)^o,$ for any $Q \in \Q_{R(B)}.$
\end{prop}

\begin{dem} Suppose that $R(B)$ is pseudo-regular and pick any $Q \in \Q_{R(B)}.$ By Lemma \ref{remark}, $(I-Q)y \in N(Q^{\#})=N(B^{\#})$ for every $y \in R(B) + R(B)^{\perpi}.$

Since $R(C) \subseteq R(B) + R(B)^{\perpi},$ we have that $B^{\#}(I-Q)C=0.$ So $X_0$ is a solution of $B^{\#}(BX-C)=0$ if and only if $B^{\#}(BX_0-QC)=0$ or $R(BX_0-QC) \subseteq R(B)^o.$ 
\end{dem}

\begin{cor} \label{Cor1013} Let $B\in CR(\HH).$ Then there exists an ImS of $BX-C=0$ for every $C \in L(\HH)$ such that $R(C) \subseteq (R(B)^o)^{\perpi}$ if and only if $R(B)$ is a pseudo-regular, nonnegative subspace of $\HH$. 
	
In this case,
$$\underset{X \in L(\HH)}{min} (BX-C)^{\#}(BX-C)=C^{\#}(I-Q)C,$$ 
for any $Q \in \Q_{R(B)}.$
\end{cor}

\begin{dem} Let $C \in L(\HH).$ Note that there exists an ImS of the equation $BX-C=0$ if and only if $R(C) \subseteq R(B) + R(B)^{\perpi}$ and $R(B)$ is nonnegative (see Theorem \ref{Teo1011}). 
	
Suppose that $R(B)$ is pseudo-regular and nonnegative. Then
$$(R(B)^o)^{\perpi}=R(B) + R(B)^{\perpi},$$ and, therefore, there exists an ImS of the equation $BX-C=0$ for every $C \in L(\HH)$ such that $R(C) \subseteq (R(B)^o)^{\perpi}.$

Conversely, suppose that there exists an ImS of the equation $BX-C=0$ for every $C \in L(\HH)$ such that $R(C) \subseteq (R(B)^o)^{\perpi}.$ Then pick a $C$ such that $R(C)=(R(B)^o)^{\perpi}=\ol{R(B) + R(B)^{\perpi}}.$ It must happen that $R(C) \subseteq R(B) + R(B)^{\perpi}$ and $R(B)$ is nonnegative. That is, $R(B)$ is to be pseudo-regular and nonnegative. In this case, let $X_0$ be an indefinite inverse of $B$ in $R(C).$ By Corollary \ref{Prop103}, $X_0$ is an ImS of $BX-C=0.$ By Proposition \ref{PropJinvPseudo}, $R(BX_0-QC) \subseteq R(B)^o.$ Then Lemma \ref{remark} with $\St=R(B)$ and the fact that $R(BX_0-QC) \subseteq R(B)^o$ yield the result.
\end{dem}

\section{Min-Max least squares problems}

In this section a min-max problem is studied for operators with not necessarly definite range. In order to pose the problem, choose a fundamental decomposition $\HH=\HH_+ \ [\dotplus] \ \HH_-$ and fix the corresponding Hilbert space $(\HH, \PI{ \ }{ \ }),$ so that, for all $x,y \in \HH$, $\PI{x }{ y }=\K{Jx}{y}$ with $J$ the signature operator associated with the decomposition.

Let $B \in CR(\HH).$ By Lemma \ref{lemaDecom}, $R(B)$ can be represented uniquely as 
\begin{equation} \label{DescompB}
R(B)=\St_+ \  [\dotplus] \ \St_-
\end{equation}
with $\St_+$ a positive closed subspace of $\HH,$ $\St_-$ a nonpositive closed subspace of $\HH$ and 
$\langle\St_+, \St_-\rangle =\{0\}.$

Consider $P_+=P_{\St_+}$ and $P_-=P_{\St_-},$ the orthogonal projections from  the Hilbert space 
$(\HH, \PI{ \ } { \ })$ onto $\St_+$ and $\St_-,$ respectively. It readily follows that $P_+ + P_- = P_{R(B)}.$ Therefore, if $B_+:=P_+B$ and $B_-:=P_-B$ then
\begin{equation} \label{DescomposicionB}
B=B_++B_-,  \quad R(B_+)=\St_+ \quad\mbox{and} \quad R(B_-)=S_-.
\end{equation}

Since $N(P_{\pm}^{\#})=\St_{\pm}^{\perpi},$ it holds that $B_+^{\#}B_-=B^{\#}P_+^{\#}B_-=0$ and $B_-^{\#}B_+=0.$ 

Observe that if $R(B)$ is regular then $P_+$ and $P_-$ are the projections given by Theorem \ref{THM2.3}. 

\begin{Def} Let  $C \in L(\HH)$. Let $B$ in $CR(\HH)$ be represented as in \eqref{DescomposicionB}. An operator $Z_0 \in L(\HH)$ is said to be an {\sl{indefinite min-max solution}} (ImMS) of $BX-C=0$ (corresponding to the decomposition given by $J$) if 
\begin{multline} \label{eqMinMax}
(BZ_0-C)^{\#}(BZ_0-C)= \\
=\underset{Y \in L(\HH) }{max} \ \ \ \underset{X \in L(\HH)}{min} (B_+X+B_-Y-C)^{\#}(B_+X+B_-Y-C).
\end{multline}
%	
%i.e., for every $x \in \HH,$
%\begin{equation} \label{eqMinMax2}
%\K{(BZ_0-C)x}{(BZ_0-C)x}=\underset{Y \in L(\HH)}{max}  \  \ \  \underset{X \in L(\HH)}{min} \K{(B_+X+B_-Y-C)x}{(B_+X+B_-Y-C)x}.
%\end{equation}
\end{Def}

The following result shows that an ImMS of $BX-C=0$ is independent of the selected fundamental decomposition  of $\HH.$ 
Along the following paragraphs, $\mathcal{C}$ denotes the cone of neutral vectors in $\HH.$

\begin{thm} \label{Teominmax} Let $C \in L(\HH)$ and $B\in CR(\HH).$ An operator $Z_0 \in L(\HH)$ is an ImMS of  $BX-C=0,$ for some (and, hence, any) fundamental decomposition of $\HH,$ if and only if $$Z_0=Z_1+Z_2$$ where $Z_1$ is an indefinite inverse of $B$ in $R(C)$ and $R(BZ_2) \subseteq \mathcal{C}.$
\end{thm}
\begin{dem} Fix a fundamental decomposition $\HH=\HH_+ \ [\dotplus] \ \HH_-,$ and consider $B=B_+ + B_-$ as in \eqref{DescomposicionB}.
Suppose that  $Z_0 \in L(\HH)$ is an ImMS of $BX-C=0$ for that decomposition. Then $Z_0$ verifies \eqref{eqMinMax}.
So, for every fixed $Y \in L(\HH),$ there exists $\underset{X \in L(\HH)}{min} (B_+X+B_-Y-C)^{\#}(B_+X+B_-Y-C).$  

From Corollary \ref{Prop103} and by using that $B_+^{\#}B_-=0,$ we get that the minimum is attained at $X_0(=X_0(Y))$ if and only if 
$$0=B_+^{\#}(B_+X_0-(C-B_-Y))=B_+^{\#}(B_+X_0-C).$$
The above says that $X_0$ is an indefinite inverse of $B_+$ in $R(C)$ and, in particular, that $X_0$ does not depend on $Y$.
Hence, for every $Y \in L(\HH),$
$$\hspace{-2.5cm}(B_+X_0+B_-Y-C)^{\#}(B_+X_0+B_-Y-C)=$$
$$= \underset{X \in L(\HH)}{min} (B_+X+B_-Y-C)^{\#}(B_+X+B_-Y-C)$$
and,  since $Z_0$ satisfies \eqref{eqMinMax}, 
% \begin{multline} \label{Eqminmax}
$$(BZ_0-C)^{\#}(BZ_0-C)
=\underset{Y \in L(\HH)}{max} (B_+X_0+B_-Y-C)^{\#}(B_+X_0+B_-Y-C).$$
%\end{multline}

By the suitable version of Corollary \ref{Prop103} and using that $B_-^{\#}B_+=0,$ we get that the maximum is attained at $Y_0 \in L(\HH)$ if and only if 
$$0=B_-^{\#}(B_-Y_0-(C-B_+X_0))=B_-^{\#}(B_-Y_0-C).$$ 
Consequently, 
\begin{equation} \label{Eqminmax2}
(BZ_0-C)^{\#}(BZ_0-C)=(B_+X_0 + B_- Y_0 - C)^{\#}(B_+X_0 + B_- Y_0 - C).
\end{equation}

Let $Z_1\in L(\HH)$  satisfy $BZ_1=B_+X_0+B_-Y_0$  as in Douglas'  Theorem, so that, according with 
\eqref{Eqminmax2}, 
\begin{equation}\label{new}
(BZ_0-C)^{\#}(BZ_0-C)=(BZ_1-C)^{\#}(BZ_1-C).
\end{equation}  
A straightforward computation gives that $B^{\#}(BZ_1-C) =0$
%$$B^{\#}(BZ_1-C)=(B_+^{\#}+B_-^{\#})(B_+X_0+B_-Y_0-C)=$$
%$$=B_+^{\#}(B_+X_0-C)+B_-^{\#}(B_-Y_0-C)=0,$$
and, in consequence, that $Z_1$ is an indefinite inverse of $B$ in $R(C)$.
Now, as $Z_1$ is an indefinite inverse of $B$ in $R(C),$ it comes that
$$
(BZ_0-C)^{\#}(BZ_0-C)=(BZ_1-C+BZ_0-BZ_1)^{\#}(BZ_1-C+BZ_0-BZ_1)=$$
$$=(BZ_1-C)^{\#}(BZ_1-C)+ (B(Z_0-Z_1))^{\#}B(Z_0-Z_1).$$
Set $Z_2:=Z_0-Z_1.$  By combining the above equation with \eqref{new} we conclude that it must hold that $(BZ_2)^{\#}BZ_2=0$ or, equivalently, that $R(BZ_2) \subseteq \mathcal{C}.$ Clearly, $Z_0=Z_1+Z_2,$ with $Z_1$  and  $Z_2$ as required.

Conversely, let $Z_1 \in L(\HH)$ be an indefinite inverse of $B$ in $R(C),$ and  $R(BZ_2) \subseteq \mc{C}.$ If $Z_0=Z_1+Z_2$ then
$$(BZ_0-C)^{\#}(BZ_0-C)=(BZ_1-C)^{\#}(BZ_1-C).$$

Write $B=B_+ + B_-$ as in \eqref{DescomposicionB}. Since $B^{\#}(BZ_1-C)=0,$ we have that $R(BZ_1-C) \subseteq N(B^{\#}) = (R(B_+) \ [\dotplus ] \ R(B_-))^{\perpi}= N(B_+^{\#}) \cap N(B_-^{\#}).$ Therefore, $B_+^{\#}(BZ_1-C)=B_-^{\#}(BZ_1-C)=0.$  
Then, as $B_-^{\#}B_+=B_+^{\#}B_-=0,$ it readily follows that, for every $X, Y \in L(\HH),$ 
$$B_+^{\#}(B_+Z_1-(C-B_-Y))=B_+^{\#}(BZ_1-C)=0$$ and $$B_-^{\#}(B_-Z_1-(C-B_+X))=B_-^{\#}(BZ_1-C)=0.$$
So, by Corollary \ref{Prop103}, we obtain that
\begin{align*}
&(BZ_0-C)^{\#}(BZ_0-C)=(BZ_1-C)^{\#}(BZ_1-C)=\\
&=(B_+Z_1+ B_-Z_1-C)^{\#}(B_+Z_1+B_-Z_1-C)\\
&=\underset{Y \in L(\HH)}{max}  (B_-Y-(C-B_+Z_1))^{\#}(B_-Y-(C-B_+Z_1))\\
&=\underset{Y \in L(\HH)}{max} \ \ \  \underset{X \in L(\HH)}{min}  (B_+X+B_-Y-C)^{\#}(B_+X+B_-Y-C).
\end{align*}
Therefore,  $Z_0$ is an ImMS of $BX-C=0.$
\end{dem}

%\bigskip
The next remark follows from the  proof of the last theorem.

\begin{obs} Let $C \in L(\HH)$ and $B\in CR(\HH)$ such that $B$ is represented as in \eqref{DescomposicionB}.  Then 
$$\underset{Y \in L(\HH)}{max}  \  \ \  \underset{X \in L(\HH) }{min} (B_+X+B_-Y-C)^{\#}(B_+X+B_-Y-C) = $$
$$ = \underset{X \in L(\HH) }{min} \  \ \ \underset{Y \in L(\HH)}{max}  (B_+X+B_-Y-C)^{\#}(B_+X+B_-Y-C).$$ 

Indeed, if $Z_0$ is an ImMS of $BX-C=0$ then, as the last theorem asserts, $Z_0=Z_1+Z_2$ where $Z_1$  is an indefinite inverse of $B$ in $R(C)$  and $R(BZ_2)\subseteq \mathcal{C}.$ In the proof of the theorem, on the other hand, we found out that, for every $X, Y \in L(\HH),$ 
	$$B_+^{\#}(B_+Z_1-(C-B_-Y))=B_+^{\#}(BZ_1-C)=0$$ and $$B_-^{\#}(B_-Z_1-(C-B_+X))=B_-^{\#}(BZ_1-C)=0.$$
A direct application of both the Corollary \ref{Prop103} and its modified version gives 
	\begin{align*}
	&(BZ_0-C)^{\#}(BZ_0-C)=(BZ_1-C)^{\#}(BZ_1-C)=\\
	&=(B_+Z_1+ B_-Z_1-C)^{\#}(B_+Z_1+B_-Z_1-C)\\
	&=\underset{X \in L(\HH)}{min}  (B_+X-(C-B_-Z_1))^{\#}(B_+X-(C-B_-Z_1))\\
	&=\underset{X \in L(\HH)}{min} \ \ \  \underset{Y \in L(\HH)}{max}  (B_+X+B_-Y-C)^{\#}(B_+X+B_-Y-C).
	\end{align*}
\end{obs}

%\vspace{0,3cm} 
\begin{cor} \label{Corollaryminmax} Let $B\in CR(\HH)$  and $C \in L(\HH).$  Then, there exists an ImMS of  $BX-C=0$  if and only if  $R(C) \subseteq R(B) + R(B)^{\perpi}.$
\end{cor}
\begin{dem} 
Suppose that $Z_0$ is an ImMS of $BX-C=0.$ Then, by Theorem \ref{Teominmax}, $Z_0=Z_1+Z_2$ where $B^{\#}(BZ_1-C)=0$ and $R(BZ_2) \subseteq \mathcal{C}.$
Therefore
$$R(C) \subseteq R(B) + R(B)^{\perpi}.$$

Conversely, if $R(C) \subseteq R(B) + R(B)^{\perpi}$ then, by Proposition \ref{PropJinvReg2}, there exists a solution of the normal equation $B^{\#}(BX-C)=0,$ say $Z_1 \in L(\HH).$  It suffices to put $Z_2=0$ and to apply Theorem \ref{Teominmax} to get that $Z_1$ is an ImMS of  $BX-C=0.$
\end{dem}

%\vspace{0,3cm} 
\begin{cor} \label{Corminmaxreg} Let $B\in CR(\HH).$  There exists an ImMS of  $BX-C=0$ for every $C \in L(\HH)$ if and only if $R(B)$ is regular.
	
	In this case, if $B$ is represented with respect  to a fixed (but arbitrary) fundamental decomposition of $\HH$  as in \eqref{DescomposicionB}, then
	$$\underset{Y \in L(\HH)}{max} \ \ \ \underset{X \in L(\HH)}{min} (B_+X+B_-Y-C)^{\#}(B_+X+B_-Y-C)=$$
	$$=C^{\#} \ [\underset{Y \in L(\HH)}{max} (B_-Y-I)^{\#}(B_-Y-I)] \ \ [ \underset{X \in L(\HH)}{min} (B_+X-I)^{\#}(B_+X-I)] \ C=$$
	$$=C^{\#}(I-Q)C,$$
	 where $Q$ is the selfadjoint projection onto $R(B).$
	
\end{cor}

\begin{dem} If $R(B)$ is regular then, for every $C \in L(\HH),$ $R(C) \subseteq R(B)  [\dotplus ] R(B)^{\perpi}$ and,  by Corollary \ref{Corollaryminmax}, there exists an ImMS of $BX-C=0.$ 

Conversely, assume that, for every $C \in L(\HH)$, there exists an ImMS of $BX-C=0.$ Set $C=I$ and apply the corollary once again to get $\HH=R(I) \subseteq R(B)+R(B)^{\perpi}$ and $R(B)$ regular.

In the case that $R(B)$ is regular, given a fundamental decomposition of $\HH,$ Ando's Theorem (Theorem \ref{THM2.3}) provides unique selfadjoint projections $Q_+, Q_- \in L(\HH)$ such that $Q=Q_+ + Q_-$ with $R(Q_+)$ uniformly positive and $R(Q_-)$ uniformly negative. Then, as we already mentioned it, the subspaces $\St_{\pm}$ in the decomposition \eqref{DescompB} of $R(B)$ and the operators $B_\pm$ in \eqref{DescomposicionB} are given by  $\St_{\pm}=R(Q_{\pm})$ and $B_{\pm}=Q_{\pm}B.$ 

Let $Z_0 \in L(\HH)$ be an ImMS of the equation $BX-C=0,$ so that, due to Theorem \ref{Teominmax}, $Z_0=Z_1+Z_2$ where $Z_1$ is an indefinite inverse in $R(C)$ and  $R(BZ_2) \subseteq \mathcal{C}.$ 
On one hand, it holds that 
$$\underset{Y \in L(\HH)}{max} \ \ \ \underset{X \in L(\HH)}{min} (B_+X+B_-Y-C)^{\#}(B_+X+B_-Y-C)$$
$$=(BZ_0-C)^{\#}(BZ_0-C)=(BZ_1-C)^{\#}(BZ_1-C)=C^{\#}(I-Q)C,$$
for $R(BZ_2) \subseteq \mathcal{C}$ and, by Corollary \ref{CorJinvC}, $BZ_1=QC.$
On the other hand, Corollary \ref{Teo101} yields 
$$\hspace{-1.5cm}C^{\#}(I-Q)C=C^{\#}(I-Q_-)(I-Q_+)C=$$
$$=C^{\#} \ [\underset{Y \in L(\HH)}{max} (B_-Y-I)^{\#}(B_-Y-I)] \ \ [ \underset{X \in L(\HH)}{min} (B_+X-I)^{\#}(B_+X-I)] \ C.$$
By merging the  above equations, the required identities are obtained and the  proof is complete.
\end{dem}

\section{The Moore-Penrose inverse in Krein spaces}

\vspace{0,3cm}  In \cite[Theorem 2.16]{XavierMary} X. Mary proved that,  given $B \in L(\HH),$ the range and nullspace of $B$ are regular subspaces of $\HH$ if and only if $B$ admits a (unique) ``Moore-Penrose inverse", in the sense that, there exists an operator $B^{\dagger} \in L(\HH)$ such that $BB^{\dagger}B=B, \ B^{\dagger}BB^{\dagger}=B^{\dagger}, \ (BB^{\dagger})^{\#}=BB^{\dagger}, \ (B^{\dagger}B)^{\#}=B^{\dagger}B.$ 

Moreover, it was proven in \cite[Corollary 2.13]{XavierMary} that if $Q$ is the selfadjoint projection onto $R(B)$ and $P$ is the selfadjoint projection onto $N(B)^{\perpi},$ then $BB^{\dagger}=Q$ and $B^{\dagger}B=P.$

In this section, we are interested in characterizing the Moore-Penrose inverse in a variational way. To this end, we consider $B\in CR(\HH)$ and $C\in L(\HH)$ and analyze the following problem: find conditions for the existence of an ImS $X_0$ of $BX-C=0$ such that $X_0^{\#}X_0 \leq Y^{\#}Y,$ for every ImS $Y$ of $BX-C=0.$ 

By Theorem \ref{Teo1011}, the equation $BX-C=0$ admits an ImS if and only if $R(C) \subseteq R(B) + R(B)^{\perpi}$ and $R(B)$ is nonnegative. In this case, if $\M_C$ is the set of ImS of $BX-C=0,$ then the above problem becomes: determine whether there exists
\begin{equation}
\underset{X \in \M_C}{min} X^{\#}X \label{eqminminC}
\end{equation}
when $\M_C \not = \emptyset.$

We only address this problem. Alternatively, symmetric problems depending on the signature of the involved subspaces can be adapted to solve them.

\begin{thm} \label{Teo411} Let $B\in CR(\HH)$ and $C \in L(\HH).$
Then there exists a solution of problem \eqref{eqminminC} if and only if $R(B)$ and $N(B^{\#}B)$ are nonnegative and $$R(C) \subseteq B(N(B^{\#}B)^{\perpi}) + R(B)^{\perpi}.$$
\end{thm}

\begin{dem} 
Suppose that there exists a solution of problem \eqref{eqminminC}.
%then there exists an ImS of $BX-C=0,$ or, equivalently, by Theorem \ref{Teo1011}, $R(B)$ is nonnegative and $R(C) \subseteq  R(B)+R(B)^{\perpi}.$
%Consider the (nonempty) set $\M_C$ of ImS of $BX-C=0.$
%Then problem \eqref{eqminminC} consists in analyzing the existence of
By Corollary  \ref{Prop103}, the set $\M_C$ can be described as 
\begin{equation*}
\M_C= \{X=X_0 + Y : Y \in L(\HH), \ R(Y) \subseteq N(B^{\#}B) \},
\end{equation*} 
where $X_0$  is any solution  of the equation $B^{\#}(BX-C)=0.$ 
	
	Therefore, problem \eqref{eqminminC} can be rephrased as: analyze the existence of
	\begin{equation}
	\underset{Z \in L(\HH)}{min} (RZ+X_0)^{\#}(RZ+X_0), \label{eq7C}
	\end{equation}
	where $R \in L(\HH)$ is such that $R(R)=N(B^{\#}B).$
	
	By Theorem \ref{Teo1011}, problem \eqref{eq7C} has a solution if and only if $N(B^{\#}B)$ is nonnegative and  $R(X_0) \subseteq  N(B^{\#}B)+N(B^{\#}B)^{\perpi}.$ Applying $B^{\#}B$ to both sides of the inclusion, we have that
	$$R(B^{\#}C)=R(B^{\#}BX_0) \subseteq B^{\#}B(N(B^{\#}B)^{\perpi}).$$ Finally, applying $(B^{\#})^{-1}$ to both sides of the inclusion, we get
	$$R(C) \subseteq B(N(B^{\#}B)^{\perpi})+N(B^{\#}).$$
	
	Conversely, suppose that $R(B)$ and $N(B^{\#}B)$ are nonnegative and \\ $R(C) \subseteq B(N(B^{\#}B)^{\perpi}) + R(B)^{\perpi}.$ Clearly, $R(C) \subseteq R(B)+R(B)^{\perpi},$ so, by Theorem \ref{Teo1011}, there exists an ImS $X_0$  of $BX-C=0,$ or equivalently, $B^{\#}(BX_0-C)=0.$ 
	On the other hand, since $N(B^{\#}B)$ is nonnegative and $R(C) \subseteq B(N(B^{\#}B)^{\perpi}) + R(B)^{\perpi}$, it holds that
	$$R(B^{\#}BX_0)=R(B^{\#}C) \subseteq B^{\#}B(N(B^{\#}B)^{\perpi}).$$ Applying $(B^{\#}B)^{-1}$ to both sides of the inclusion, it comes that
	$$R(X_0) \subseteq N(B^{\#}B) + N(B^{\#}B)^{\perpi}.$$ Therefore, by Theorem \ref{Teo1011}, there exists a solution of \eqref{eq7C} and hence, there exists a solution of problem  \eqref{eqminminC}.
\end{dem}

It follows from the last theorem that, if $B\in CR(\HH)$ and $C \in L(\HH)$ are such that $R(B)$ and $N(B^{\#}B)$ are nonnegative, then there exists a solution of problem \eqref{eqminminC} if and only if $\M_C  \not = \emptyset,$ and for every $X_0 \in \M_C,$ $R(X_0) \subseteq N(B^{\#}B) + N(B^{\#}B)^{\perpi}.$ Moreover: 

\begin{lema} \label{Prop1050} Let $B\in CR(\HH)$ and $C \in L(\HH)$ such that $R(B)$ and $N(B^{\#}B)$ are nonnegative and $R(C) \subseteq B(N(B^{\#}B)^{\perpi}) + R(B)^{\perpi}.$
	Then $X_1$ is a solution of \eqref{eqminminC} if and only if $B^{\#}(BX_1-C)=0$ and $R(X_1) \subseteq N(B^{\#}B)^{\perpi}.$
\end{lema}

\begin{dem} Recall that $X_1$ is a solution of problem \eqref{eqminminC} if and only if  $X_1=RZ_1 + X_0,$ with $R \in L(\HH)$ such that $R(R) = N(B^{\#}B),$  $X_0$ a solution of $B^{\#}(BX-C)=0$ and $Z_1$  a solution of \eqref{eq7C}.
	
	Since $N(B^{\#}B)$ is nonnegative, by Theorem \ref{Teo1011} and Corollary \ref{Cor1012}, $Z_1 \in L(\HH)$ is a solution of \eqref{eq7C} if and only if $Z_1$ is such that
	$$R^{\#}(RZ_1+X_0)=0$$ that is, $$R^{\#}X_1=0$$ or, equivalently, 
	$$R(X_1) \subseteq N(R^{\#}) = R(R)^{\perpi} = N(B^{\#}B)^{\perpi}.$$
\end{dem}

As a corollary of Theorem \ref{Teo411}, we have the following result.

\begin{prop} \label{Teo41} Let $B\in CR(\HH).$ Then the following assertions are equivalent:
\begin{enumerate}
	\item [i)] There exists a solution of problem \eqref{eqminminC} for $C=I,$
	\item [ii)] $R(B)$ and $N(B)$ are uniformly positive,
	\item [iii)] there exists the Moore-Penrose inverse of $B,$ $B^{\dagger},$ and $R(B)$ and $N(B)$ are nonnegative.
\end{enumerate}

\end{prop}

\begin{dem} $i) \Leftrightarrow ii):$ If there exists a solution of problem \eqref{eqminminC} for $C=I$, by Theorem \ref{Teo411}, $R(B)$ and $N(B^{\#}B)$ are nonnegative, and 
\begin{equation} \label{eqregular}	
\HH \subseteq B(N(B^{\#}B)^{\perpi}) + R(B)^{\perpi}.
\end{equation}
Then, clearly, $R(B)$ is regular and nonnegative, i.e., $R(B)$ is uniformly positive.
Since $R(B)$ is regular then $R(B^{\#})=R(B^{\#}B)$ or, equivalently, $N(B)=N(B^{\#}B).$ 
Applying $B^{\#}$ to both sides of \eqref{eqregular}, we have that
$R(B^{\#}) \subseteq B^{\#}B(N(B^{\#}B)^{\perpi}).$ Then, 
$$R(B^{\#}B)=R(B^{\#}) \subseteq B^{\#}B(N(B^{\#}B)^{\perpi}) \subseteq R(B^{\#}B).$$ Therefore
$$R(B^{\#}B)=B^{\#}B(N(B^{\#}B)^{\perpi}),$$ and so,
$$\HH=N(B^{\#}B)+N(B^{\#}B)^{\perpi}=N(B)+N(B)^{\perpi}.$$ Thus, $N(B)$ is regular and nonnegative and therefore uniformly positive.

Conversely, if $R(B)$ and $N(B)$ are uniformly positive, then 
$$\HH=R(B)+R(B)^{\perpi}=N(B^{\#}B)+N(B^{\#}B)^{\perpi},$$ where we used the fact that $N(B)=N(B^{\#}B)$ since $R(B)$ is regular. Applying $B^{\#}B$ to both sides of the second equality, we get that $R(B^{\#}B)=B^{\#}B(N(B^{\#}B)^{\perpi}).$ Then, applying $(B^{\#})^{-1}$ to the left and right sides of the last equality, the inclusion
$R(B) + R(B)^{\perpi} \subseteq B(N(B^{\#}B)^{\perpi}) + R(B)^{\perpi}$ holds. Therefore we get
$$\HH= B(N(B^{\#}B)^{\perpi}) + R(B)^{\perpi},$$ and, by Theorem \ref{Teo411}, there exists a solution of problem \eqref{eqminminC} for $C=I.$

$ii) \Leftrightarrow iii):$ See \cite[Theorem 2.6]{XavierMary}.
\end{dem}

%\begin{prop} \label{Prop1050} Let $(\HH, \K{ \ }{ \ })$ be a Krein space and $B\in L(\HH)$ such that $R(B)$ and $N(B)$ are both uniformly positive.
%Then $X_1$ is solution of problem \eqref{eqminmin} if and only if $B^{\#}(BX_1-I)=0$ and $R(X_1) \subseteq N(B)^{\perpi}.$
%\end{prop}
%
%\begin{dem} Recall that $X_1$ is a solution of problem \eqref{eqminmin} if and only if  $X_1=PZ_1 + X_0,$ with $P \in L(\HH)$ such that $R(P) = N(B),$  $X_0$ a solution of $BX=Q,$ with $Q$ the selfadjoint projection onto $R(B),$ and $Z_1$  a solution of \eqref{eq7}.
%	
%Since $N(B)$ is uniformly positive, by Theorem \ref{Teo1011} and Corollary \ref{Cor1012}, $Z_1 \in L(\HH)$ is a solution of \eqref{eq7} if and only if $Z_1$ is such that
%$$P^{\#}(PZ_1+X_0)=0$$ that is, $$P^{\#}X_1=0,$$ or, equivalently 
%$$R(X_1) \subseteq N(P^{\#}) = R(P)^{\perpi} = N(B)^{\perpi}.$$
%\end{dem}

\begin{thm} \label{Teo66} Let $B\in CR(\HH)$ and suppose that $N(B)$ and $R(B)$ are uniformly positive.  Then, the Moore-Penrose inverse of $B,$ $B^{\dagger},$ is the unique ImS $X_0$ of $BX-I=0$ such that $X_0^{\#}X_0 \leq Y^{\#}Y,$ for every ImS $Y$ of $BX-I=0.$ 
\end{thm}

\begin{dem}  Since $N(B)$ and $R(B)$ are regular, the Moore-Penrose inverse of $B,$ $B^{\dagger},$ exists. Consider $Q$ the selfadjoint projection onto $R(B)$ and $P$ the selfadjoint projection onto $N(B)^{\perpi},$ then
$$B^{\#}(BB^{\dagger}-I)=B^{\#}(Q-I)=0.$$ 
On the other hand, $$R(B^{\dagger})=R(B^{\dagger}BB^{\dagger})\subseteq R(B^{\dagger}B) =  R(P) = N(B)^{\perpi}.$$ Hence, by Lemma \ref{Prop1050}, $B^{\dagger}$ is a solution of problem \eqref{eqminminC}, with $C=I.$
	
Let $X_1 \in L(\HH)$ be any other solution of problem \eqref{eqminminC}, with $C=I.$ By Lemma \ref{Prop1050}, $X_1$ is an ImS of $BX-I=0$ and $R(X_1) \subseteq N(B)^{\perpi}.$ Then, by Corollary \ref{CorJinvC}, $BX_1=Q=BB^{\dagger}.$  Therefore, 
$$X_1=PX_1=B^{\dagger}BX_1=B^{\dagger}BB^{\dagger}=B^{\dagger}.$$ 

\end{dem}

The next remark follows from the proofs of Proposition \ref{Teo41} and Theorem \ref{Teo66}.
\begin{obs}  Let $B\in CR(\HH)$ and $C \in L(\HH),$ and suppose that $N(B)$ and $R(B)$ are uniformly positive.  
Then problem \eqref{eqminminC} admits a unique solution, namely, $B^{\dagger}C.$

%In fact, if $N(B)$ and $R(B)$ are uniformly positive, then $$\HH=R(B)+R(B)^{\perpi}=N(B^{\#}B)+N(B^{\#}B)^{\perpi},$$ where we used the fact that $N(B)=N(B^{\#}B)$ since $R(B)$ is regular. Then, by similar arguments as in Proposition \ref{Teo41}, we have that 
%$$R(C) \subseteq R(B)+R(B)^{\perpi} \subseteq B(N(B^{\#}B)^{\perpi}) + R(B)^{\perpi}.$$ Therefore, by Theorem \ref{Teo411}, there exists a solution of problem \eqref{eqminminC}. 
%
%Moreover, since $N(B)$ and $R(B)$ are regular, the Moore-Penrose inverse of $B,$ $B^{\dagger},$ exists. Let $X_1=B^{\dagger}C,$ then proceeding as in the last theorem, it follows that $X_1$ is the unique solution of problem \eqref{eqminminC}. 

\end{obs}

\subsection{The Moore-Penrose inverse: the pseudo-regular case}

In \cite[Proposition 5.1]{GiribetIndefinite}, a family of generalized inverses of a closed range operator with pseudo-regular range and nullspace was given. In this case, the associated projections turn out to be normal. In this section, we prove the equivalence between the existence of this family of generalized inverses and the pseudo-regularity of the range and nullspace of an operator $B \in CR(\HH).$ 
We also give a more general expression for these generalized inverses and we characterize them in a variational way as we did in the last section with the Moore-Penrose inverse.

Given $B\in CR(\HH),$ recall that $\tilde{B}$ is a $\{1,2\}$-inverse of $B$ if $\tilde{B}$ is a solution of the system
$$
BXB=B, \ XBX=X.
$$
If $(\HH, \PI{\cdot}{\cdot})$ is a Hilbert space, every $B\in CR(\HH)$ admits a $\{1,2\}$-inverse, see \cite[Theorem 3.1]{AriasGI}. Then, using any of the underlying Hilbert structures, the same is true in the Krein space $\HH.$
Observe that, if $\tilde{B}$ is a $\{1,2\}$-inverse of $B,$ then $B\tilde{B}$ is a projection onto $R(B)$ and $\tilde{B}B$ is a projection with $N(\tilde{B}B)=N(B).$

\begin{prop} \label{proppseudo} Let $B\in CR(\HH).$  Then, there exists a solution of the system 
	\begin{equation} \label{EqPseudo}
\left \{
\begin{array}{ccc} 
BXB=B,\ XBX=X,\\
(BX)^{\#}(BX)=(BX)(BX)^{\#},\\
(XB)^{\#}(XB)=(XB)(XB)^{\#},
\end{array} 
\right.
\end{equation}
if and only if $R(B)$ and $N(B)$ are pseudo-regular subspaces of $\HH.$

In this case, $D \in L(\HH)$ is a solution of \eqref{EqPseudo} if and only if there exist $Q \in \Q_{R(B)}$ and $P \in \Q_{N(B)}$ such that 
\begin{equation} \label{eqD}
D=(I-P) \tilde{B} Q,
\end{equation} where $\tilde{B}$ is any $\{1,2\}$-inverse of $B.$
\end{prop}

\begin{dem} Suppose that $R(B)$ and $N(B)$ are pseudo-regular subspaces. Let $Q \in \Q_{R(B)}$ and $P \in \Q_{N(B)}.$ Let $\tilde{B}$ be any $\{1,2\}$-inverse of $B.$ 
Let $D$ be defined as in \eqref{eqD}. From $BP=0$ it follows immediately that
$$BD=Q, \mbox{ and } DB=I-P.$$ So the last two equations of the system are satisfied.
Also, 
$$BDB=QB=B \mbox{ and } DBD=(I-P)D=D.$$

Conversely, suppose that \eqref{EqPseudo} admits a solution $D.$ Let $Q=BD$ and $P=I-DB,$ then $P$ and $Q$ are normal projections in $L(\HH).$ Moreover, 
$R(Q)=R(BD) \subseteq R(B).$ On the other hand, $R(Q)=R(BD) \supseteq R(BDB)=R(B).$ Therefore, $R(Q)=R(B)$ and $R(B)$ is pseudo-regular. Also, $N(B) \subseteq N(DB) = N(P) \subseteq N(BDB) =N(B).$ So that $N(B)=N(P)$ and then $N(B)$ is pseudo-regular.

In this case, we have already proven that if $D$ is as in \eqref{eqD}, then $D$ is a solution of \eqref{EqPseudo}. Conversely, suppose that $D \in L(\HH)$ is a solution of \eqref{EqPseudo}. Note that $Q:=BD \in \Q_{R(B)}$ and $P:=I-DB \in \Q_{N(B)}.$ Let  $\tilde{B}$ be any $\{1,2\}$-inverse of $B.$ It is straightforward to check that $D$ satisfies 
$$ (I-P)\tilde{B}Q=D.$$
\end{dem}

\begin{prop} \label{Cor411} Let $B\in CR(\HH),$ such that $R(B)$ is pseudo-regular. Then there exists a solution of problem \eqref{eqminminC} for every $C\in L(\HH)$ such that $R(C) \subseteq \ol{B(N(B^{\#}B)^{\perpi}) + R(B)^{\perpi}}$ if and only if $N(B^{\#}B)$ is pseudo-regular and \break $N(B^{\#}B)$ and $R(B)$ are nonnegative.
\end{prop}

\begin{dem} Suppose that $R(B)$ and $N(B^{\#}B)$ are nonnegative and pseudoregular. Then, by \cite[Lemma 3.4]{GiribetIndefinite}, $R(B^{\#}B)$ is closed. Since $N(B^{\#}B)$ is pseudo-regular, \cite[Corollary 2.5]{Izumino} gives that  $R(B^{\#}BB^{\#}B)$ is closed too.
	
Suppose that $R(C) \subseteq \ol{B(N(B^{\#}B)^{\perpi}) + R(B)^{\perpi}}.$ Then 
$$R(B^{\#}C) \subseteq B^{\#} [(\ol{B(R(B^{\#}B)) + R(B)^{\perpi}})]$$
$$\subseteq B^{\#} [(R(B) \cap R(BB^{\#}B)^{\perpi})^{\perpi}]\subseteq [B^{-1}(R(B) \cap R(BB^{\#}B)^{\perpi})]^{\perpi} =$$
$$ = [B^{-1}(R(BB^{\#}B)^{\perpi})]^{\perpi}=[B^{\#}R(BB^{\#}B)]^{\perpi  \ \perpi}=R(B^{\#}BB^{\#}B),$$ where we used the fact that $B^{\#}(\St^{\perpi}) \subseteq (B^{-1}(\St))^{\perpi},$ for any closed subspace $\St \subseteq \HH$ and that $R(B^{\#}BB^{\#}B)$ is closed. 
Then, applying $(B^{\#})^{-1}$ to both sides of the inclusion, we have that
$$R(C) \subseteq (B^{\#})^{-1}(R(B^{\#}BB^{\#}B))=N(B^{\#})+ B(N(B^{\#}B)^{\perpi}).$$ Whence, by Theorem \ref{Teo411}, problem \eqref{eqminminC} admits a solution.

%Conversely, suppose that there exists a solution of problem \eqref{eqminminC}, for every $C\in L(\HH)$ such that $R(C) \subseteq \ol{B(N(B^{\#}B)^{\perpi}) + R(B)^{\perpi}},$ then there exists a solution of problem \eqref{eqminminC}, for every $C\in L(\HH)$ such that $R(C) \subseteq \ol{R(B) + R(B)^{\perpi}},$ therefore there exists an ImS of $BX-C=0$, for every $C\in L(\HH)$ such that $R(C) \subseteq \ol{R(B) + R(B)^{\perpi}},$ then by Corollary \ref{Cor1013}, $R(B)$ is regular and nonnegative.

Conversely, suppose that there exists a solution of problem \eqref{eqminminC} for every $C\in L(\HH)$ such that $R(C) \subseteq \ol{B(N(B^{\#}B)^{\perpi}) + R(B)^{\perpi}}.$  Then pick $C$ such that $$R(C) = \ol{B(N(B^{\#}B)^{\perpi}) + R(B)^{\perpi}}.$$ By Theorem \ref{Teo411}, we have that $N(B^{\#}B)$ and $R(B)$ are nonnegative and $R(C)=\ol{B(N(B^{\#}B)^{\perpi}) + R(B)^{\perpi}} \subseteq B(N(B^{\#}B)^{\perpi}) + R(B)^{\perpi}.$ Then the subspace $B(N(B^{\#}B)^{\perpi}) + R(B)^{\perpi}$ is closed. Hence,
$$B^{-1}(B(N(B^{\#}B)^{\perpi}) + R(B)^{\perpi})=N(B)+N(B^{\#}B)^{\perpi}+N(B^{\#}B)=$$
$$=N(B^{\#}B)^{\perpi}+N(B^{\#}B)$$ is closed, so that $N(B^{\#}B)$ is pseudo-regular.
\end{dem}

The next result is a corollary of Proposition \ref{proppseudo}. We will use it in the proof of Theorem \ref{TeoPseudo} in order to characterized the solutions of \eqref{eqminminC} in term of pseudo-inverses when $R(B)$ and $N(B^{\#}B)$ are pseudo-regular.

\begin{lema} \label{CorPseudo}  Let $B\in CR(\HH)$ such that $R(B)$ is a pseudo-regular subspace of $\HH.$ 
	Given $Q \in \Q_{R(B)},$ let $B'=Q^{\#}B.$ Then there exists a solution of the system 
	\begin{equation} \label{EqPseudo2}
	\left \{
	\begin{array}{ccc} 
	B'XB'=B',\ XB'X=X,\\
	B'X=Q^{\#}Q,\\
	(XB')^{\#}(XB')=(XB')(XB')^{\#},
	\end{array} 
	\right.
	\end{equation}
	if and only if $N(B^{\#}B)$ is a pseudo-regular subspace of $\HH.$
	
	In this case, $D \in L(\HH)$ is a solution of \eqref{EqPseudo2} if and only if there exists $P \in \Q_{N(B^{\#}B)}$ such that $$D=(I-P) \tilde{B'} Q^{\#}Q,$$ where $\tilde{B'}$ is any $\{1,2\}$-inverse of $B'.$
\end{lema}

\begin{dem} Note that 
	$$R(B')=R(Q^{\#}Q) \mbox{ and } N(B')=N(B^{\#}B).$$
	Then apply Proposition \ref{proppseudo} to $B'.$
\end{dem}

\begin{thm} \label{TeoPseudo} Let $B\in CR(\HH)$ and $C \in L(\HH).$ If $R(B)$ and $N(B^{\#}B)$ are nonnegative pseudo-regular subspaces and  $R(C) \subseteq \ol{B(N(B^{\#}B)^{\perpi}) + R(B)^{\perpi}},$ set $X_1=DC,$ where $D\in L(\HH)$ is a solution of  \eqref{EqPseudo2}, then $X_1$ is a solution of problem \eqref{eqminminC}.
\end{thm}

\begin{dem} By the proof of Proposition \ref{Cor411}, the set  $B(N(B^{\#}B)^{\perpi}) + R(B)^{\perpi}$ is closed.
	
Given a solution $D$ of \eqref{EqPseudo2}, consider $P \in \Q_{N(B^{\#}B)},$ $Q \in \Q_{R(B)}$ and any $\{1,2\}$-inverse $\tilde{B'}$ of $B'$ such that 
$$D=(I-P) \tilde{B'} Q^{\#}Q.$$ Observe that 
$$Q^{\#}BDC=B'DC=Q^{\#}QC=Q^{\#}C,$$ where we used the fact that $R(C) \subseteq  R(B)+ R(B)^{\perpi}$ and Lemma \ref{remark}.
Then $R(BDC-C) \subseteq N(Q^{\#})=N(B^{\#})$ or, equivalently,
$$B^{\#}(BDC-C)=0.$$ Then, by Proposition \ref{PropJinvPseudo}, $DC$ is an ImS of $BX-C=0.$ 

On the other hand, $$R(B^{\#}BDC)=R(B^{\#}C) \subseteq B^{\#}(B(N(B^{\#}B)^{\perpi})),$$ so, by applying $(B^{\#}B)^{-1}$ to both sides of the inclusion, we have that
$$R(DC) \subseteq N(B^{\#}B) + N(B^{\#}B)^{\perpi} = N(P^{\#}(I-P)).$$
Then $$P^{\#}(I-P)DC=P^{\#}DC=0.$$ Thus $R(DC) \subseteq N(B^{\#}B)^{\perpi}$ and, by Lemma \ref{Prop1050}, $X_1=DC$ is a solution of problem \eqref{eqminminC}.
\end{dem}

\begin{obs} Under the same assumptions of the last theorem, by Proposition \ref{Cor411}, there exists a solution of problem \eqref{eqminminC}.  Furthermore, if $R(C) \not \subseteq R(B)^{\perpi},$ a converse of Theorem \ref{TeoPseudo} holds:  if $X_1$ is a solution of problem \eqref{eqminminC} then $X_1=DC,$ where $D\in L(\HH)$ is a solution of \eqref{EqPseudo2}.
	
In fact, let $X_1$ be a solution of problem \eqref{eqminminC}, then by similar arguments as those in \cite[Theorem 3.5]{GiribetIndefinite}, there exists $P' \in \Q_{N(B^{\#}B)}$ such that
$$X_1=(I-P')X_0,$$ where $X_0$ is an ImS of $BX-C=0.$

Let $Q \in \Q_{R(B)}$ and $\tilde{B'}$ be any $\{1,2\}$-inverse of $B'.$ Set
$$D=(I-P') \tilde{B'} Q^{\#}Q.$$ Then, by Lemma \ref{CorPseudo}, we have that $D$ is a solution of \eqref{EqPseudo2}. Then, proceeding as in the proof of the last theorem, we get that $DC$ is an ImS of $BX-C=0.$
Then, by Corollary \ref{Prop103}, $X_0=DC+Y,$ with $R(Y) \subseteq N(B^{\#}B).$ Hence, $$X_1=(I-P')X_0=(I-P')DC=DC.$$
\end{obs}

\section*{Acknowledgements}
Maximiliano Contino was supported by CONICET PIP 0168.  Alejandra Maestripieri was partially supported by CONICET PIP 0168. The work of Stefania Marcantognini was done during her stay at  the Instituto Argentino de Matem\'atica with an appointment funded by the CONICET. She is greatly grateful to the institute for the hospitality and to the CONICET for financing her post.

\end{document}